\documentclass{gtart_h}  

\def\ifplaintex{\expandafter\ifx\csname documentclass\endcsname\relax}

\def\ifplaintex{\expandafter\ifx\csname documentclass\endcsname\relax}

\ifplaintex 
\else
\fi

\expandafter\ifx\csname epsfbox\endcsname\relax\input epsf\fi

\def\gt{{\mathsurround=0pt\it $\cal G\mskip-2mu$eometry \&\ 
$\cal T\!\!$opology}}        %  journal title in recommended style

\def\gtp{{\mathsurround=0pt\it $\cal G\mskip-2mu$eometry \&\ 
$\cal T\!\!$opology $\cal P\!$ublications}}  % GT publications

\def\lognumber#1{\def\thelognumber{#1}}
\def\volumenumber#1{\def\thevolumenumber{#1}}
\def\papernumber#1{\def\thepapernumber{#1}}
\def\volumeyear#1{\def\thevolumeyear{#1}}

\def\pagenumbers#1#2{\def\startpage{#1}\def\finishpage{#2}}
\def\published#1{\def\publishdate{#1}}
\def\proposed#1{\def\theproposer{#1}}
\def\seconded#1{\def\theseconders{#1}}
\def\received#1{\def\receiveddate{#1}}

\def\accepted#1{\def\accepteddate{#1}}

\let\\\par\let\thelognumber\relax
\let\thevolumenumber\relax\let\thepapernumber\relax
\let\thevolumeyear\relax\let\thesamplenumber\relax\let\startpage\relax
\let\finishpage\relax\let\publishdate\relax\let\receiveddate\relax
\let\reviseddate\relax\let\accepteddate\relax\let\theasciititle\relax
\let\theasciiauthors\relax
\let\theasciiabstract\relax
\let\theasciiemail\relax\let\theshortauthors\relax\let\theshorttitle\relax

\long\def\maketitlep{   % start of definition of \maketitlep

\count0=\startpage

\gt\hfill      %   Journal title (top left) 
\hbox to 77pt{\vbox to 0pt{\vglue -15pt\epsfbox{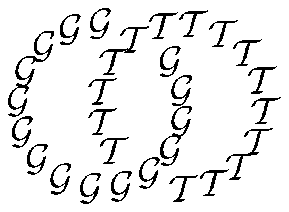}\vss}\hss}
\break
{\small\ifx\thesamplenumber\relax % sample?  
Volume \else Sample
\fi\thevolumenumber\ (\thevolumeyear)
\startpage--\finishpage\nl
Published: \publishdate}
\vglue 0.5truein plus 0.4fil minus 0.1truein

{\parskip=0pt\leftskip 0pt plus 1fil\def\\{\par\smallskip}{\ifplaintex\large
\else\Large\fi\bf\thetitle}\par\medskip}   

\vglue 0pt plus 0.1fil 

{\parskip=0pt\leftskip 0pt plus 1fil\def\\{\par}{\sc\theauthors}
\par\medskip}

\vglue 0pt plus 0.1fil 

{\small\parskip=0pt\let\newline\\
{\leftskip 0pt plus 1fil\def\\{\par}{\sl\theaddress}\par}
\expandafter\ifx\theemail\relax    % email address?
\relax\else\vglue 5pt plus 0.02fil minus 2pt\def\\{\stdspace{\rm 
and}\stdspace} 
\cl{Email:\stdspace\tt\theemail}\fi
\ifx\theurl\relax                  % URL given?
\relax\else\vglue 5pt plus 0.02fil minus 2pt\def\\{\stdspace{\rm 
and}\stdspace}
\cl{URL:\stdspace\tt\theurl}\fi\par}

\vglue 7pt plus 0.3fil minus 3pt

{\bf Abstract}
\vglue 5pt plus 0.1fil minus 2pt

\theabstract

\vglue 7pt plus 0.3fil minus 3pt

{\bf AMS Classification numbers}\quad Primary:\quad \theprimaryclass

Secondary:\quad \thesecondaryclass

\vglue 5pt plus 0.3fil minus 2pt

{\bf Keywords:}\quad \thekeywords

\vglue 10pt plus 0.5fil minus 5pt

{\small  Proposed: \theproposer\hfill Received: \receiveddate\nl
Seconded: \theseconders\hfill 
\ifx\reviseddate\relax                         % paper revised?
Accepted: \accepteddate                        % no
\else
Revised: \reviseddate                          % yes
\fi}
\eject
}       %  end of definition of \maketitlep

\font\phead=cmsl9 scaled 950
\font=cmsl9 scaled 1050
\font\pnum=cmbx10 scaled 913
\font\lnum=cmbx10 
\font\pfoot=cmsl9 scaled 950
\font=cmsl9 scaled 1050
\ifplaintex
\headline{\vbox to 0pt{\vskip -4.5mm\line{\small\phead\ifnum
\count0=\startpage ISSN 1364-0380 (on line)
1465-3060 (printed) \hfill {\pnum\folio}\else\ifodd\count0\def\\{ }% 
\ifx\theshorttitle\relax\thetitle\else\theshorttitle\fi\hfill{\pnum\folio}
\else\def\\{ and }{\pnum\folio}\hfill\ifx\theshortauthors\relax\theauthors
\else\theshortauthors\fi\fi\fi}\vss}}
\footline{\vbox to 0pt{\vglue 0mm\line{\small\pfoot\ifnum\count0=\startpage
\copyright\ \gtp\hfill\else
\gt, Volume \thevolumenumber\ (\thevolumeyear)\hfill\fi}\vss
}}
\else
\makeatletter
\def\@oddhead{{\small\lhead\ifnum\count0=\startpage ISSN 1364-0380 (on line)
1465-3060 (printed) \hfill {\lnum\number\count0}\else\ifodd\count0
\def\\{ }\ifx\theshorttitle\relax \thetitle \else\theshorttitle\fi\hfill
{\lnum\number\count0}\else\def\\{ and }{\lnum\number\count0}
\hfill\ifx\theshortauthors\relax 
\theauthors\else\theshortauthors\fi\fi\fi}}\def\@evenhead{\@oddhead}
\def\@oddfoot{\small\lfoot\ifnum\count0=\startpage\copyright\ \gtp\hfill\else
\gt, Volume \thevolumenumber\ (\thevolumeyear)\hfill\fi}
\def\@evenfoot{\@oddfoot}
\makeatother
\fi

\newwrite\gtoutfile
\long\gdef\makeheadfile{  %%% start of definition of \makeheadfile
{\def\\{, }\def\s{ }
\immediate\openout\gtoutfile head.xxx
\immediate\write\gtoutfile{Proxy-for: \ifx\theasciiauthors\relax
\theauthors\else\theasciiauthors\fi\s<\ifx\theasciiemail\relax\theemail\else\theasciiemail\fi>}
\immediate\write\gtoutfile{\noexpand\\}
\immediate\write\gtoutfile{Authors: \ifx\theasciiauthors\relax
\theauthors\else\theasciiauthors\fi}
{\def\\{ }\immediate\write\gtoutfile{Title: \ifx\theasciititle\relax
\thetitle\else\theasciititle\fi}}
\immediate\write\gtoutfile{Subj-class: GT or SG or MG etc}
\immediate\write\gtoutfile{MSC-class: \theprimaryclass\ifx\thesecondaryclass\relax\else, \thesecondaryclass\fi}
\immediate\write\gtoutfile{Journal-ref: Geom. Topol. \thevolumenumber
(\thevolumeyear) \startpage-\finishpage}
\immediate\write\gtoutfile{Comments: Published by Geometry and Topology at}
\immediate\write\gtoutfile{\s\s http://www.maths.warwick.ac.uk/gt/GTVol\thevolumenumber/paper\thepapernumber.abs.html}
\immediate\write\gtoutfile{\noexpand\\}
\immediate\write\gtoutfile{}
\ifx\theasciiabstract\relax
\immediate\write\gtoutfile{\theabstract}\else
\immediate\write\gtoutfile{\theasciiabstract}\fi
\immediate\write\gtoutfile{}
\immediate\write\gtoutfile{\noexpand\\}
\immediate\write\gtoutfile{}
\immediate\closeout\gtoutfile}}  %%% end of definition of \makeheadfile

\def\maketitlepage{\maketitlep\makeheadfile}
\let\maketitle\maketitlepage

\lognumber{434}
\received{17 March 2004}
\volumenumber{9}\papernumber{4}\volumeyear{2005}
\pagenumbers{163}{178}   
\published{6 January 2005\nl Corrected: 7 March 2005}
\accepted{27 December 2004}

\proposed{Martin Bridson}
\seconded{Steven Ferry, Benson Farb}

\usepackage{amsmath,amssymb}
 
\let\Bbb\mathbb
\def\eagerbreak{\par \ifdim\lastskip<\medskipamount \removelastskip
\penalty-175 \vskip\medskipamount\fi}
\newenvironment{XX}{\eagerbreak\refstepcounter{XXX}%
\noindent\par\noindent{\bf \theXXX\qua}\sl}{\smallbreak}
\newenvironment{XY}{\eagerbreak\refstepcounter{XXX}%
\noindent\par\noindent{\bf \theXXX\qua}}{\smallbreak}
\newenvironment{Enumerate}{\begin{list}{\bf(\roman{enum})}%
{\usecounter{enum}}}{\end{list}}
\newenvironment{Thm}{\begin{XX}{\bf Theorem\qua}\ignorespaces}{\end{XX}}
\newenvironment{Prop}{\begin{XX}{\bf Proposition\qua}\ignorespaces}{\end{XX}}
\newenvironment{Defn}{\begin{XY}{\bf Definition\qua}\ignorespaces}{\end{XY}}
\newenvironment{Cor}{\begin{XX}{\bf Corollary\qua}\ignorespaces}{\end{XX}}
\newenvironment{Lemma}{\begin{XX}{\bf Lemma\qua}\ignorespaces}{\end{XX}}

\newenvironment{Pf}{\proof}{\endproof}
\newenvironment{Rmk}{\begin{XY}{\bf Remark\qua}\ignorespaces}{\end{XY}}

\newenvironment{Example}{\begin{XY}{\bf Example\qua}\ignorespaces}{\end{XY}}
{\end{Enumerate}}
{\eqno{\rm (\theXXX)} $$}  
{\setcounter{equation}{\value{XXX}} \begin{eqnarray}}%
{\end{eqnarray} \setcounter{XXX}{\value{equation}}}
\newcommand{\Twoline}[2]{$$\displaylines{\quad{#1}\hfill\cr\hfill{#2}\quad\cr}%
$$}
\newcommand{\Opencone}{{\mathcal O}}
\newcommand{\Supp}{\operatorname{\rm Supp}}
\newcommand{\Prob}{\operatorname{\rm Prob}}

\newcommand{\R}{\mbox{$\Bbb R$}}
 
\newcounter{XXX}[section]
\renewcommand{\theXXX}{\thesection.\arabic{XXX}}

\renewcommand{\phi}{\varphi}
\renewcommand{\epsilon}{\varepsilon}
\newcounter{enum}

\newcommand{\itm}{}

\newcounter{ritmctr}
\newenvironment{ritem}{\begin{itemize}\setcounter{ritmctr}{0}%
 \renewcommand{\itm}{\addtocounter{ritmctr}{1}\item[\rm(\roman{ritmctr})]}}%
{\end{itemize}}
\newcounter{aitmctr}
\newenvironment{aitem}{\begin{itemize}\setcounter{aitmctr}{0}%
 \renewcommand{\itm}{\addtocounter{aitmctr}{1}\item[\rm(\alph{aitmctr})]}}%
{\end{itemize}}

\begin{document}
\bibliographystyle{gtart}
 
\title{Warped cones and property A}

\author{John Roe}
\address{Department of Mathematics, Penn State University\\University
Park PA 16802, USA}
\email{roe@math.psu.edu}
\begin{abstract}
We describe a construction (the `warped cone construction') which produces
examples of coarse spaces with large groups of translations.  We show that by
this construction we can obtain many examples of coarse spaces which do not have
property A or which are not uniformly embeddable into Hilbert space.
\end{abstract}
\primaryclass{53C20}
\secondaryclass{43A07, 53C12, 20F69}
\keywords{Coarse geometry, amenable action, property A, warped cone}
\maketitlepage

The purpose of this paper is to explore some aspects of the geometry of those
coarse spaces which may be created from compact foliated manifolds, or from
discrete group actions on compact metric spaces, by way of the \emph{warped
cone construction}.  This construction was introduced in~\cite{JR17} (see
also~\cite{JR18}) and some of its properties were predicted there. Since those
papers were written, workers in coarse geometry such as Yu~\cite{Yu6} have
drawn attention to the importance of coarse (or uniform) \emph{embeddability}
in Hilbert space, and the related `coarse amenability' \emph{property A}.  We
shall show that the warped cone construction is sufficiently flexible to create
natural new examples both of spaces which have, and of spaces which do not
have, property A\null.  In a sequel to this paper~\cite{HR6} we will show that the
class of warped cones also includes counterexamples to the coarse Baum--Connes
conjecture.

The first section of this paper reviews in detail the construction of warped
cones (the article~\cite{JR17} contains no proofs). In the second section we
briefly review and generalize the definition of property A; for the purposes of
this paper it is necessary to extend the definition from the context of
uniformly discrete spaces, in which it is presented in the literature, to that
of more general coarse spaces.  In the third section we show that the warped
cone over an amenable action is a space with property A\null.  The fourth and final
section contains a partial converse:   if the warped
cone over a particular kind of action has property A, then the action itself is
amenable.  This allows us to construct warped cones that do not have property A,
or that are not uniformly embeddable into Hilbert space.
There is a striking parallel between the theory of warped cones and the theory
of the `box spaces' of~\cite[Chapter 11]{JR27}; this will be brought out in the
main text.   

It is a pleasure to thank Eric Guentner and Nigel Higson for illuminating
discussions.

\section{Definition and basic properties of warped cones}

Let $(X,d)$ be a proper metric space.  We are going to consider $X$ as a coarse
space (see \cite{HPR} for the definition), so that only the `large scale'
properties of the metric on $X$ will be relevant to us.  Recall that a metric
space is \emph{proper} if closed bounded subsets are compact.

\begin{Example} We will be particularly concerned throughout this paper with the
case $X=\Opencone Y$, the \emph{open cone} on a compact space $Y$.  We will
assume that $Y$ is a `nice' space such as a smooth manifold or a finite
simplicial complex.  The open cone $X=\Opencone Y$ may be defined in two
equivalent ways, which we state for $Y$ a manifold:
\begin{aitem}
\itm Embed $Y$ smoothly in a high-dimensional sphere $S^{N-1}$, and let $X$ be
the union of all the rays through the origin in $\R^N$ that meet the embedded
copy of $Y$; equip $X$ with the metric (distance function) induced from $\R^N$.
\itm Let $X= Y\times [1,\infty) $ as a manifold\footnote{Here we have---for
simplicity---omitted a neighborhood of the cone point.  Our interest is in the
other end of the cone, in the large scale behavior of the geometry, so that the
treatment of the singularity does not matter.}, equipped with the
Riemannian metric $t^{2}g_Y + g_{\mathbb R}$, where $g_Y$ is a Riemannian
metric on $Y$, $g_{\mathbb R}$ is the standard Riemannian metric on $\R$, and $t$ is
the coordinate on $\R$.
\end{aitem}
(There are obvious `piecewise linear' versions of both constructions in case $Y$
is a finite complex.) It is a straightforward exercise to show that up to coarse
equivalence these constructions are independent of the choices involved and,
moreover, that both constructions yield coarsely equivalent spaces $X$.
\end{Example}

\begin{Defn} A map $f\colon X\to X$ is a \emph{translation} if it is a bijection
and if
the supremum
$ \sup_{x\in X} d(x,f(x) )$ 
is finite. \end{Defn}

The translations of a coarse space form a group which is closely related to the
coarse structure of the space.

\begin{Example} Let $\Gamma$ be a discrete group equipped with a \emph{right}
invariant metric.  Then the \emph{left} action of $\Gamma$ on itself is an
action by translations. \end{Example}

Suppose now that $X$ is a coarse space and that $\Gamma$ is a group
acting\footnote{Actions of a group on a space will be taken on the left.} on
$X$ (the example of particular interest to us is $X=\Opencone Y$, $\Gamma$ a
group of diffeomorphisms of $Y$.)  It will typically not be the case that
$\Gamma$ acts by translations on $X$ and we want to modify the metric on $X$ in
such a way as to ensure that the action of $\Gamma$ becomes an action by
translations.  We will describe this process as \emph{warping} $X$ by $\Gamma$.

\begin{Lemma} Let $\mathcal F$ be a family of metrics on a set $X$.  If the supremum
\[ d(x,x')= \sup \{ \delta(x,x'): \delta\in{\mathcal F} \} \]
is finite for each pair $(x,x')\in X\times X$, then it defines a metric on $X$.
\end{Lemma}

The proof is elementary. 

\begin{Defn} Let $(X,d)$ be a proper metric space and let $\Gamma$ be a group
acting by homeomorphisms on $X$, provided with a finite generating set $S$.  The \emph{warped
metric} $d_\Gamma$ on $X$ is the greatest metric that satisfies the inequalities
\[ d_\Gamma(x,x') \le d(x,x'), \qquad d_\Gamma(x, sx) \le 1 \ \forall\ s\in S.
\]
\end{Defn}

The existence and uniqueness of the warped metric follow easily from the Lemma
above---the warped metric is the supremum of the family $\mathcal F$ of
metrics satisfying the inequalities that appear in its definition.  More
constructively, we can say:

\begin{Prop}\label{chainprop}Let $x,y\in X$.  For $\gamma\in\Gamma$, let $|\gamma|$ denote the
word length of $\gamma$ relative to the generating set $S$. 
 The warped distance from $x$ to $y$ is the
infimum of all sums
\[ \sum d(\gamma_ix_i, x_{i+1}) + |\gamma_i| \]
taken over all finite sequences $x=x_0, x_1, \ldots, x_N = y$ in $X$ and
$\gamma_0, \ldots, \gamma_{N-1}$ in $\Gamma$.  Moreover, if $d(x,y)\le k$, then
one can find a finite sequence as above with $N=k+1$ such that the infimum is
attained. \end{Prop}

\begin{Pf} The expression in the proposition defines a metric, call it $\delta$,
that belongs to the family of metrics $\mathcal F$ mentioned above, and thus
$\delta\le d_\Gamma$ by definition.  On the other hand, $d_\Gamma\le\delta$ by
the triangle inequality.  Thus the two metrics are equal.  This proves the first
part of the proposition.

For the second part, note that in order for the sum above to be $\le k+1$, no
more than $k+1$ group elements $\gamma_i$ can differ from the identity.  If say
$\gamma_i$ is equal to the identity, then $\gamma_i$ and $x_i$ may be omitted
from the chain and by the triangle inequality the sum will only decrease.  Thus
every chain may be shortened to one having $N\le k+1$.  The set of such chains
beginning at $x$ and ending at $y$ is compact (by the properness of the metric
space $X$), so the infimum is attained. \end{Pf}

\begin{Prop} The warped metric is a proper metric.  The coarse structure induced
by the warped metric does not depend on the choice of generating set $S$ for
$\Gamma$, nor on the choice of metric $d$ within the coarse structure of $X$.
\label{warpmprop}\end{Prop}

\begin{Pf} For a subset $K$ of $X$, define $N_k(K)$ to be  
\[ N(K) = \{ x\in X: d(x,\gamma y)\le 1 \mbox{\  for some $y\in K$, $|\gamma |
\le 1$}\}; \]
that is, $N_k(K)$ is the union of the closed $k$--neighborhoods of translates of $K$
by generators of the group, together with the closed $k$--neighborhood of $K$
itself.  Since $X$ is proper and $\Gamma$ acts by homeomorphisms, the operation
$N_k$ preserves compactness: if $K$ is compact, then so is $N_k(K)$.  It follows
from Proposition~\ref{chainprop} that the warped $k$--ball around a point $p$ is
contained in the compact set $(N_k)^{k+1}\{p\}$.  This proves properness.

Consider the abstract coarse structure $\mathcal{S}$ on $X$ generated by the controlled sets
for the original metric on $X$ together with the sets $\{(x,\gamma x):
\gamma\in\Gamma\}$.  This structure is countably generated and hence metrizable
\cite[Theorem 2.55]{JR27}.  The generators for $\mathcal{S}$ are all controlled for
the warped metric, so every $\mathcal{S}$--controlled set is controlled for the warped
metric.  On the other hand, by Proposition~\ref{chainprop} every set controlled
for the warped metric is $\mathcal{S}$--controlled.  Thus the warped metric describes
the coarse structure $\mathcal{S}$.  Since the definition of $\mathcal{S}$ does not
involve the choice of metric on $X$ nor the choice of generating set for
$\Gamma$, any two warped metrics are coarsely equivalent. 
\end{Pf}

We can now make the key definition of this paper.

\begin{Defn} Let $Y$ be a smooth compact manifold, or a finite simplicial
complex, and let $\Gamma$ be a finitely generated group acting by
homeomorphisms on $Y$. The \emph{warped cone} $\Opencone_\Gamma Y$ is the coarse
space obtained by warping $\Opencone Y$ along the induced $\Gamma$--action.
\end{Defn}

Typically the metric structure of $\Opencone_\Gamma Y$ is difficult to envisage.
 Each cross-section $Y_t$ is a copy of $Y$ with the metric modified in the
 following manner: first the ambient metric of $Y$ is multiplied by $t$, and
 then `shortcuts' or `wormholes' are introduced along the generators of the
 group action, each of which `costs' one unit of distance to travel through.  It
 is clear that as $t$ increases, the distortion introduced by the warping
 becomes progressively more severe.   Nevertheless, the $Y_t$ do increase in
 size:  

\begin{Lemma} With notation as above,
  the diameter of $Y_t$ tends to infinity as $t\to\infty$, unless $Y$ is 
  a finite set
  and $\Gamma$ acts transitively on it. Moreover, if $\Gamma$ acts by
  diffeomorphisms on the compact manifold $Y$, the diameter of $Y_t$ grows at
  least as fast as $\log t$. \end{Lemma}

\begin{Pf} The proof of the first part is left to the reader.  For the 
second part, use elementary estimates of volume: the volume, in some
fixed Riemannian metric on $Y$, of the warped $d$--ball in $Y_t$ about
some point, is of order at most the number of group elements of word
length $\le d$ times the maximum volume of a set of the form $\gamma
(B)$, where $B$ is a $d/t$-ball in $Y$ and $\gamma$ is a word in
$\Gamma$ of length at most $d$.  This volume is at most $C e^{kd}
(d/t)^n$, where $C$ and $k$ are constants and $n$ is the dimension of
$Y$.  If this quantity is less than the volume of $Y$, then $Y_t$ has
diameter $\ge d$; the logarithmic growth follows. \end{Pf}

\begin{Prop} Suppose that $\Gamma$ acts on $Y$ through  Lipschitz
homeomorphisms.  Then the
 space $X=\Opencone_\Gamma Y$ has bounded geometry. \end{Prop}

\begin{Pf}  (I am grateful for Graham Niblo and Nick Wright for pointing out
an error in the proof given in the first version of this paper.)  Let
$X$ be a metric space and let $\epsilon>0$.  By definition, the
\emph{$\epsilon$-capacity} of a subset $A$ of $X$ is the maximum
cardinality of an $\epsilon$-separated subset of $A$.  We shall say
that $X$ has \emph{bounded geometry} if there is an $\epsilon>0$ such
that for every $r>0$ there is a $s>0$ such that each $r$-ball in $X$
has $\epsilon$-capacity at most $s$.  An open cone $X=\Opencone Y$ on
a manifold or finite complex has bounded geometry, as one sees most
easily by regarding it as a metric subspace of a Euclidean space.

Suppose that $X$ has bounded geometry.  We will say for the purposes
of this proof that a subset of $X$ has \emph{size $\le (p,q)$} if it
is the union of at most $p$ subsets of $X$ each of which has diameter
at most $q$.  Clearly, for each $(p,q)$ there is $C=C(p,q)$ such that
any subset of size $\le (p,q)$ has $\epsilon$-capacity at most $C$.
 
Now let $X=\Opencone Y$, where $\Gamma$ acts on $Y$ (and therefore on
$X$) through Lipschitz homeomorphisms.  From its definition, there
exist constants $C$ and $c$ such that if $S$ has size $\le (p,q)$
(relative to the unwarped metric on $X$) , then $ N_k(S)$ has size
$\le (Cp, cq+2k)$; $C$ is the number of generators for the group, and
$c$ is the maximum Lipschitz constant for the action of a generator.
But it was shown in the proof of Proposition~\ref{warpmprop} that the
$d_\Gamma$-ball of radius $k$ around any point $p$ is contained in
$(N_k)^{k+1}\{p\}$.  It therefore follows by induction that the
$\epsilon$-capacity (in the $d$-metric) of any such ball is bounded.
However, $d_\Gamma\le d$ and therefore the capacity of this ball in
the $d_\Gamma$-metric is at most equal to its capacity in the
$d$-metric.  Bounded geometry now follows.
\end{Pf}  

\begin{Rmk} In \cite{JR17} we defined warped cones not for 
\emph{group actions} but for
\emph{foliations}.  We review this construction briefly. 
Let $V$ be a compact smooth manifold equipped with a foliation $F$.  Choose a
`normal bundle' $N$ to the   foliation, that is a complementary subbundle to $TF$
in $TM$.  Choose Euclidean metrics $g_N$ in $N$ and $g_F$ in $TF$.  The
\emph{foliated warped cone} $\Opencone_F V$ is the manifold $V\times
[0,\infty)/V\times \{0\}$
equipped with the distance function induced for $t\ge 1$ by the Riemannian metric $g_{\mathbb R} +
g_F + t^2g_N$.  As in our earlier discussion, the structure near $t=0$ does not
much matter from the point of view of coarse geometry; we imagine that the space
$\Opencone_F V$ has been equipped with some arbitrary path metric in this
region.
One sees without difficulty that the coarse structure of $\Opencone_F V$ does not depend on the
choice of normal bundle $N$, nor does it depend on the choice of metrics $g_F$
and $g_N$ (compare Proposition~\ref{warpmprop}). \end{Rmk}

Foliations are often constructed by `suspending' group actions.   For foliations
constructed in this way, there is a relationship between our two notions of
warped cone.   Let $\Gamma$ be a group acting freely by diffeomorphisms on a manifold
$M$, and suppose that $\Gamma$ has a classifying space $B\Gamma$ which is itself
a compact manifold.  Form the balanced product $V=M\times_\Gamma E\Gamma$.  This
is a compact manifold. The foliation of $M\times E\Gamma$ by copies of $E\Gamma$
descends to a foliation $F$ of $V$ whose leaves are also copies of
$E\Gamma$ (since we are assuming that the action is free).  The manifold $M$ may
be embedded in $V$ as a transversal to the foliation.

\begin{Lemma}\label{equivlemma}In the above situation, the embedding of $M$ as a transversal in
$V$ induces a coarse equivalence $\Opencone_\Gamma(M)\to \Opencone_F(V)$.
\end{Lemma}

We leave the proof to the reader. 

\section{Property A for general spaces}

  For a locally compact,
 $\sigma$--compact topological space $X$ (e.g. a bounded geometry metric space),
let $\Prob(X)$ denote the space of Radon probability measures on $X$.
It is the state space of the $C^*$--algebra $C_0(X)$ and as such it carries two
topologies: the norm topology and the weak-$*$ topology.  We will use both in
the definition below.

\begin{Defn}\label{Adef}Let $X$ be a bounded geometry proper metric space. We say that $X$
has \emph{property A} if there exists a sequence of weak-$*$ continuous maps
$f_n\colon X\to\Prob(X)$ such that
\begin{ritem}
\itm for each $n$ there is an $r$ such that, for each $x$, the measure $f_n(x)$
is supported within $B(x;r)$, and
\itm for each $s>0$, as $n\to\infty$, 
\[ \sup_{d(x,y)<s} \|f_n(x)-f_n(y)\| \to 0. \]
\end{ritem}
For brevity we will refer to {\rm (i)} above as the property of being
\emph{uniformly localized}. \end{Defn}

If $X$ is a \emph{discrete} metric space of bounded geometry, then this is Yu's
definition as reformulated in \cite[Lemma 3.5]{HR5}.  Moreover, there is a
simple reduction to the discrete case.  A uniformly discrete subset $Z\subseteq
X$ will be called a \emph{lattice} in $X$ if there is $R>0$ such that
$X=\bigcup_{z\in Z} B(z;R)$.  The inclusion $Z\to X$ is then a coarse
equivalence.

\begin{Lemma} A (bounded geometry proper metric) space $X$ has property A if and
only if some (and hence every) lattice in $X$ has property A.
\end{Lemma}

\begin{Pf} Let $Z$ be a lattice in $X$ and let $\{\phi_z\}$ be a partition of
unity on $X$ subordinate to a cover of $X$ by balls $B(z;R)$, $z\in Z$.

Suppose that $X$ has property $A$ and let maps $f_n\colon X\to\Prob(X)$ be 
as in the definition~\ref{Adef}.  Define maps $g_n\colon Z\to\Prob(Z)$ by
\[ g_n(z) = \sum_{w\in Z} \left( \int \phi_w df_n(z) \right) \delta_w, \]
where $\delta_w$ denotes the Dirac mass at $w$.  This is a  
probability measure, uniformly localized near $z$, and we have
\Twoline{\|g_n(z)-g_n(z')\| = \sum \left| \int \phi_w d[f_n(z)-f_n(z')]
\right|}{ 
\le \sum\int \phi_w d|f_n(z)-f_n(z')| = \|f_n(z)-f_n(z')\|.}
This shows that $Z$ has property A.

Conversely, suppose that $Z$ has property A and let $g_n\colon Z\to\Prob(Z)$ be
a defining sequence of maps.  We may consider probability measures on Z as
probability measures on X (linear combinations of Dirac measures).  Set
\[ f_n(x) = \sum_z \phi_z(x) g_n(z) \]
where $g_z$ is considered as a probability measure on $X$.  Then $f_n(x)$ is a
probability measure, uniformly localized near $x$, depending continuously on
$x\in X$, and we have
\[ \|f_n(x)-f_n(x')\| \le \sup_{z,z'\in B(x;R)\cup B(x';R)} \|g_n(z)-g_n(z')\|
\]
since both $f_n(x)$ and $f_n(x')$ are convex combinations of the measures
$g_n(z)$, for $z\in B(x;R)\cup B(x';R)$. This shows that $X$ has property
A. 
\end{Pf}

It is known in the discrete case that there is an alternative characterization
of property A in terms of kernels of positive type.  We will also need the
continuous analogue of this characterization.  Recall that a continuous function
$k\colon X\times X\to\R$ is a \emph{(continuous) kernel of positive type} if 
\[ \sum_{i,j=1^n} \lambda_i\lambda_j k(x_i,x_j) \ge 0 \]
for all $n=1,2,\ldots$ and all
 $n$--tuples $\lambda_1,\ldots,\lambda_n$ in $\R$ and $x_1,\ldots,x_n$ in
$X$.   

A kernel $k$ has \emph{controlled support} if there is $R>0$ such that
$k(x,x')=0$ whenever $d(x,x')>R$.

\begin{Prop} The (bounded geometry proper metric) space $X$ has property A if
and only if there is a sequence $\{k_n\}$ of   continuous positive definite
kernels    such that 
\begin{aitem}
\itm $|k_n(x,x')|\le 1$ for all $x,x'$;
\itm each $k_n$ has controlled support;
\itm 
$ k_n(x,x')\to 1$ as $n\to\infty$, 
\emph{uniformly} on each controlled set $\{(x,x'): d(x,x')<r\}$.\end{aitem}
\label{positiveA} \end{Prop}

\begin{Pf}  We shall take this proposition for granted in the uniformly discrete
case (see~\cite{JR27} or \cite{Tu1} for the details) so our task is simply to show
that $X$ has the three-part property (temporarily call it `P') of the Proposition if and only if some lattice $Z$ in
$X$ has it.

If $X$ has  P, then we may simply restrict the positive definite
kernels $k_n$ to $Z$ to show that $Z$ has P also.

Conversely, suppose that $Z$ has P.  Let $\{\phi_z\}$ be a partition of unity
subordinate to some cover $B(z;R)$, as in the proof of the previous proposition.
If $\ell_n$ is a sequence of positive definite kernels exhibiting property P for
$Z$, then one checks without difficulty that
\[ k_n(x,x') = \sum_{z,z'} \ell_n(z,z') \phi_z(x)\phi_{z'}(x')  \]
is a sequence of such kernels exhibiting property P for $X$.
\end{Pf}

It will be helpful to consider an apparent weakening of
property $A$.

\begin{Defn} Let $X$ be a bounded geometry proper metric space. We say that $X$
has \emph{property A at infinity} if there exists a sequence of weak-$*$ continuous maps
$f_n\colon X\to\Prob(X)$ such that
\begin{ritem}
\itm for each $n$ there is an $r$ such that, for each $x$, the measure $f_n(x)$
is supported within $B(x;r)$, and
\itm for each $s>0$, as $n\to\infty$, 
\[ \inf_{\strut\mbox{\scriptsize\rm $K\subseteq X$ compact}}\quad
\sup_{ d(x,y)<s,\ x,y\notin K} \|f_n(x)-f_n(y)\| \to 0. \]
\end{ritem}
\end{Defn}

Clearly, if a space has property A, it has property A at infinity.  It is
useful to observe that the converse is true also.

\begin{Lemma} If a space $X$ has property A at infinity, then it has property
A.
\end{Lemma}
\begin{Pf} Let $\{f_n\}$ be a sequence of maps provided by the definition of
property A at infinity.  We are going to construct 
a sequence $g_m$ of weak-$*$
continuous maps $X\to\Prob(X)$ which are uniformly localized and satisfy
\[ \sup_{d(x,y)\le m} \|g_m(x)-g_m(y)\| \le m^{-1}. \]
To do so, fix $m$ and then choose $n$ so large that for some compact  
$K\subseteq X$,
\[ \sup_{ d(x,y)<m,\ x,y\notin K} \|f_n(x)-f_n(y)\| \le (4m)^{-1}; \]
we can do this by definition of property $A$ at infinity.  Fix a
function $\phi\colon X\to [0,1]$ such that $\phi(x)=1$ if $x\in K$, $\phi(x)=0$
if $d(x,K)>5m^2$, and $\phi$ is Lipschitz with Lipschitz constant less than
$(2m)^{-2}$; for instance, $\phi$ might be a suitable function of $d(x,K)$. 
 Let
$\delta_K$ be the
Dirac measure supported at some (arbitrarily chosen) point of $K$. Now
put
\[ g_m(x) = \phi(x)\delta_K + (1-\phi(x))f_n(x). \]
It is clear that $g_m$ is a uniformly localized map from $X$ to $\Prob(X)$.
Moreover if $d(x,y)<m$ then 
\[ \|g_m(x) - g_m(y)\| \le \frac{1}{4m} + \frac{1}{4m} + \|f_n(x)-f_n(y)\| \le
\frac{3}{4m} \]
and so we are done.
\end{Pf}

\section{Warped cones over amenable actions}

\begin{Prop} Let $(X,d)$ be a bounded geometry coarse space that has property
A, and let $\Gamma$ be an amenable group acting on $X$ (by coarse maps).   
  Then the warped space $(X,d_\Gamma)$ has property
A also. \label{warpA}\end{Prop}

\begin{Cor} Let $Y$ be a compact manifold (or finite simplicial complex) and let
$\Gamma$ be an amenable group acting by Lipschitz homeomorphisms
on $Y$.  Then the warped cone
$\Opencone_\Gamma(Y)$ has property A. \label{cbcamen}\end{Cor}

\begin{Pf}
Let $f_n\colon X\to\Prob(X)$ be the functions that are provided by the
definition of property A\null.  Let $\mu_n$ be a F\o lner sequence in
$\Prob(\Gamma)$, by which we mean that $\{\mu_n\}$ is a sequence of finitely
supported probability measures on $\Gamma$ such that, for each fixed
$\gamma\in\Gamma$, the difference $\mu_n - \gamma\cdot\mu_n$ tends to zero in
norm as $n\to\infty$.\footnote{The action of $\Gamma$ on $\Prob(\Gamma)$ is
defined by
$(\gamma\cdot\mu)\{\gamma'\} = \mu\{\gamma' \gamma^{-1}\}$.}

For fixed $m>0$ choose $n_1$ so large that $\|\mu_{n_1} - \gamma\cdot\mu_{n_1}\|
< (2m)^{-2}$ for all $\gamma\in\Gamma$ having $|\gamma|\le m$.  Let $R$ be defined
by
\[ R = \sup \{ d(\gamma x, \gamma y) : d(x,y)\le m,\ \gamma\in\Supp(\mu_{n_1})
\} ; \]
$R$ is finite because the supremum is taken over a finite set of coarse maps.
Finally choose $n_2$ so large that 
$ \|f_{n_2}(x) - f_{n_2}(y) \| < (2m)^{-2}$
whenever $d(x,y)\le R$. 
Define a weak-$*$ continuous function $g_m$ from $X$ 
to $\Prob(X)$ by convolving $\mu_{n_1}$ with
$f_{n_2}$:
\[ g_m(x) = \sum_{\gamma\in\Gamma} \mu_{n_1}\{\gamma\} f_{n_2}(\gamma x). \] 
It is clear
from the construction of the warped metric $d_\Gamma$ that   $g_m$ is
$d_\Gamma$--uniformly localized; in fact $g_m(x)$ is supported within the $d_\Gamma$--ball
around $x$ of radius  
\[ \inf \{ r: \Supp(f_{n_2}(x))\subseteq B_d(x;r)\,\forall x\} + \max \{ |\gamma| :
\gamma\in\Supp(\mu_{n_1}) \} , \]
the first term being finite by the definition of property A and the second being
finite by our assumption that $\mu_{n_1}$ is finitely supported.
 We will show that
  $\|g_m(x)-g_m(y)\|\le 1/m$   for $d_\Gamma(x,y)<m$; this
will complete the proof that the maps $\{g_m\}$ are a sequence realizing property A
for the $d_\Gamma$--metric.

Because of Proposition~\ref{chainprop} it suffices to demonstrate 
that $\|g_m(x)-g_m(y)\|\le 1/(2m^2)$ in two cases:
 first in the case that $d(x,y)\le m$, and
second in the case that $y = \gamma' x$ with
 $\gamma'\in\Gamma$, $|\gamma'|\le m$.  In the first case
we may write
\[ \|g_m(x) - g_m(y)\| \le \sum_\gamma \mu_{n_1}\{\gamma\} \|f_{n_2}(\gamma x) -
f_{n_2}(\gamma y)\|. \]
By hypothesis we have $d(\gamma x, \gamma y)\le R$ for each term appearing on
the right and therefore each term $\|f_{n_2}(\gamma x) -
f_{n_2}(\gamma y)\|$ is bounded by $(2m)^{-2}$.   Since
$\sum_\gamma \mu_n\{\gamma\}=1$ the desired result follows.  In the second case write,
by a simple change of variable,
\[\| g_m(x) - g_m(\gamma' x) \| \le
 \sum_\gamma \Bigl| \mu_{n_1}\{\gamma\} - \mu_{n_2}\{\gamma
(\gamma')^{-1}\}\Bigr| \| f_{n_2}(\gamma x)\|. \]
By hypothesis,
 \[ \sum_\gamma \Bigl| \mu_{n_1}\{\gamma\} - \mu_{n_2}\{\gamma
(\gamma')^{-1}\}\Bigr| \le (2m)^{-2}. \]
 This gives the result in the second case also.
\end{Pf}

Let $\Gamma$ be a (discrete) group that acts on a compact  Hausdorff space $Y$. 
Recall that one says that the action of $\Gamma$ is \emph{amenable} if there
exists a sequence of weak-$*$ continuous maps $\mu_n\colon Y \to \Prob(\Gamma)$ such that, for each
$\gamma\in\Gamma$,
\[ \sup_{y\in Y} \|\gamma\cdot \mu_n(y) - \mu_n(\gamma y)\| \to 0 \]
as $n\to\infty$.  Thus a group $\Gamma$ is amenable if and only if its action on
a point is amenable.  We may and shall assume that for each $n$ there is some
finite subset of $\Gamma$ on which all the measures $\mu_n(y)$ are supported.

\begin{Rmk} It is clear that the definitions of property A and of amenable action are
closely related. In fact, the main result of~\cite{HR5} is that a group $\Gamma$
has property A if and only if its natural action on its Stone--Cech
compactification is amenable. \end{Rmk}

We want to generalize Corollary~\ref{cbcamen} to the case of amenable
\emph{actions}. 

\begin{Thm} Suppose that $Y$ is a compact manifold, or a finite simplicial
complex, and that the group $\Gamma$ acts on $Y$ amenably by Lipschitz
homeomorphisms.
 Then the warped
cone $X=\Opencone_\Gamma(Y)$ has property A.    \label{amenacttheorem}\end{Thm}

\begin{Pf} It will suffice to prove property A at infinity.  The proof is an
elaboration of the proof of Proposition~\ref{warpA}.
Let $f_n\colon X\to\Prob(X)$ be the functions that are provided by the
definition of property A\null.  Let $\mu_n$ be a F\o lner sequence of maps
$Y\to \Prob(\Gamma)$, as in the definition of amenable action, and extend them
radially to maps $X\to\Prob(\Gamma)$ (because we are only concerned with
property A at infinity we do not need to worry about what happened near the cone
point).
As before, define functions $g$ on $X$ by `convolving' $\mu$ with
$f$:
\[ g_m(x) = \sum_{\gamma\in\Gamma} \mu_{n_1}(x)\{\gamma\} f_{n_2}(\gamma x),\] 
where $n_1$ and $n_2$ are chosen by the same   procedure as in the earlier
proof. 
As before 
these are $d_\Gamma$--uniformly localized \mbox{weak-$*$} continuous maps from 
$X$ to $\Prob(X)$.  To complete the proof of property A at infinity we must
demonstrate the `approximate invariance' condition
that for fixed $s$, \[
\inf_{\strut\mbox{\scriptsize\rm $K\subseteq X$ compact}}\quad
\sup_{ d(x,y)<s,\ x,y\notin K} \|g_m(x)-g_m(y)\| \to 0 \]
as $n\to\infty$.

Because of Proposition~\ref{chainprop} it suffices to demonstrate the desired
convergence separately in two cases: first in the case that $d(x,y)\le s$, and
second in the case 
that $y = \gamma' x$ with $\gamma'\in\Gamma$, $|\gamma'|\le s$.
  In the first case
we may write
\Twoline{ \|g_m(x) - g_m(y)\| \le \sum_\gamma \mu_{n_1}(x)\{\gamma\} \|f_{n_2}(\gamma x) -
f_{n_2}(\gamma y)\| }{+ \sum_\gamma |\mu_{n_1}(x)\{\gamma\} - \mu_{n_1}(y)\{\gamma\}|
\|f_{n_2}(\gamma y)\|.}
The first sum on the right-hand side tends to zero by the same argument as
before.  As for the second sum, notice that because of the weak-$*$ continuity
of   $\mu_{n_1}$ when considered as a map on $Y$, together with the geometry of the
$d$--metric on the open cone, there is for any $\epsilon>0$ a compact subset $K$
of $X$ such that $|\mu_{n_1}(x)\{\gamma\} - \mu_{n_1}(y)\{\gamma\}| < \epsilon$ for
$x,y\notin K$, $d(x,y)<s$.  Thus the second sum may be made arbitrarily small by
restricting $x,y$ to lie outside a compact set $K$.  This completes the proof in
the first case.

 In the second case $y=\gamma' x$ write,
by a simple change of variable,
\[ g_m(x) - g_m(\gamma' x) = \sum_\gamma \Bigl( \mu_{n_1}(x)\{\gamma\} -
\mu_{n_1}(\gamma' x)\{\gamma
(\gamma')^{-1}\}\Bigr) f_{n_2}(\gamma x). \]
By definition of a F\o lner sequence the term in parentheses tends to 0 in norm,
uniformly for $|\gamma' |\le s$.  This gives the result in the second case also.
\end{Pf}

\section{Non property A examples}

In this section we will consider examples of the following sort.  Let $G$ be a
compact Lie group, and let $\Gamma\le G$ be a dense finitely generated subgroup, which
we consider as a discrete group acting on $G$ by left translation.  We study the
coarse geometry of the warped cone $\Opencone_\Gamma(G)$.

\begin{Prop} Suppose that the warped cone $\Opencone_\Gamma(G)$, as defined
above, has property A\null.  Then $\Gamma$ is amenable.
\end{Prop}

\begin{Pf}
Since $X=\Opencone_\Gamma(G)$ has property A, there is a sequence $\{k_n\}$ of
continuous positive type kernels on $X$ satisfying the conditions of
Proposition~\ref{positiveA}.  By averaging with respect to the Haar measure on
$G$, we may assume that these kernels are invariant under left multiplication by
$G$;
in particular, they are invariant under the (left translation) action of
$\Gamma$.  Thus for each $t>0$ there is a function $h_{n,r}$ on $\Gamma$ of
positive type such that
\[ k_n((r,g),(r,\gamma g)) = h_{n,r}(\gamma) \]
for every $g\in G$.

The functions $h_{n,r}$ are all bounded in absolute value by 1 and therefore we
can find a function $h_n$ on $\Gamma$ of positive type such that
\[ h_{n,r_k}(\gamma) \to h_n(\gamma) \]
for a suitable sequence $r_k\to\infty$.

Fix $\gamma$.  For all sufficiently large $r$, the warped distance from $(r,g)$
to $(r,\gamma g)$ is precisely the word length of $\gamma$.  Since $k_n$ has
controlled support, it follows that $h_n(\gamma)=0$ whenever the word length of
$\gamma$ is sufficiently large.  In other words, $h_n$ has compact support. 
Furthermore, since $k_n\to 1$ uniformly on controlled sets as $n\to\infty$, we
have $h_n\to 1$ pointwise.  Thus $\Gamma$ admits a  sequence of positive definite
functions of compact support tending pointwise to 1, and this is a well-known
characterization of amenability. \end{Pf}

\begin{Rmk} From this result and Theorem~\ref{amenacttheorem} it follows that if
$\Gamma$ acts amenably on $G$, then $\Gamma$ itself is amenable.  This is easy
to prove directly. \end{Rmk}

\begin{Rmk} It is well known that $G=SO(3)$ contains a free subgroup $\Gamma$ 
on two generators.  Thus $\Opencone_\Gamma(G)$ is a coarse space that does not
have property A. \end{Rmk}

By a similar device we may prove:

\begin{Prop} Suppose that the warped cone $\Opencone_\Gamma(G)$, as defined
above, is uniformly embeddable in Hilbert space.
  Then $\Gamma$ has the Haagerup property~\cite{CCJJV}.  In particular, if $\Gamma$ has
  property T, then $\Opencone_\Gamma(G)$ cannot be uniformly embedded in Hilbert
  space.
\end{Prop} 

Examples of dense property T subgroups of compact Lie groups can be found for
instance with $G=SO(5)$; see~\cite{Lubotzky}.  In the sequel we will show that
some of the resulting warped cones are counterexamples to the coarse Baum--Connes
conjecture.

\begin{Pf} We use the same idea as the preceding proof.  To say that
$X=\Opencone_\Gamma(G)$ admits a uniform embedding into Hilbert space is to say
that $X$ admits a continuous negative type kernel $k$ which is \emph{effective},
that is, there are unbounded increasing functions
$\rho_1,\rho_2\colon\R^+\to\R^+$ such that
\[ \rho_1(d(x,x'))\le k(x,x') \le \rho_2(d(x,x')). \]
Averaging over $G$, we may assume that $k$ is $G$--invariant.  For a fixed
$\gamma$, 
\[ h_r(\gamma) = k((r,x),(r,\gamma x)) \]
depends only on $\gamma$ and is bounded by $\rho_2(|\gamma|)$; so we can find a
sequence $r_k\to\infty$ such that the limit 
\[ h(\gamma) = \lim_k h_{r_k}(\gamma) \]
exists.  It is a negative type function and satisfies $h(\gamma)\ge
\rho_1(|\gamma|)$, so it is proper.   Thus $\Gamma$ admits a proper negative
type function, which is the Haagerup property. \end{Pf}

Compare these results with Propositions 11.26 and 11.39 in~\cite{JR27}.

\end{document}